\documentclass[11pt]{amsart}
\usepackage{amsfonts,amssymb,stmaryrd,amscd,amsmath,latexsym,amsbsy}
\newtheorem{thm}{Theorem}

\newtheorem{prop}{Proposition}

\newtheorem{defi}{Definition}
\newtheorem{lem}{Lemma}

\newtheorem*{sol}{Solution}

\newcommand{\solu}[1]{\begin{sol}{\bf (\ref{#1})}}

\newcommand{\CC}{\mathbb C}

\title{Finite dimensional representations of the double affine Hecke
algebra of rank $1$}
\author{A.Oblomkov}\author{E.Stoica}
\address{Department of Mathematics, Princeton University, Fine Hall, 
Washington Road, Princeton, NJ 08544, USA.}
\address{Department of Mathematics, MIT, 77, Massachusetts Ave., Cambridge,
MA 02139, USA.}
\email{oblomkov@math.princeton.edu}\email{immanuel@math.mit.edu}
\begin{document}
\begin{abstract}
We classify the finite dimensional irreducible representations of
the double affine Hecke algebra (DAHA) of type $C^{\vee}C_1$ in the
case when $q$ is not a root of unity.
\end{abstract}

\maketitle

\section{Introduction}

The double affine Hecke algebras (DAHA) were introduced by Cherednik
at the beginning of nineties \cite{C1},\cite{C2}. Roughly speaking,
DAHA is a nontrivial extension of the affine Hecke algebra
corresponding to an affine root system  by the group algebra of the
coroot lattice. Later, it was discovered that these algebras were
connected to many areas of mathematics. For instance, DAHAs lie at
the core of the proof of the Macdonald conjectures, given by
Cherednik \cite{C2}. Subsequently, Etingof and Ginzburg \cite{EG}
revealed deep connections between the theory of the rational
degenerations of DAHA and  noncommutave geometry.

In this paper, we study the finite dimensional complex
representations of the DAHA of rank $1$, namely the algebra
$H=H(k_0,k_1,u_0,u_1;q)$ corresponding to the root system
${CC_1}^{\vee}$. This algebra was defined by Sahi \cite{Sa} (see
also \cite{NS}). As it was shown in \cite{O}, it is the most general
DAHA of rank $1$.

We restrict our attention to the case when $q$ is not a root of
unity. The case when $q=1$ and the parameters $k_0,k_1,u_0,u_1$ are
generic was studied in \cite{O}. In the special case when
$k_0=u_0=u_1=1$, the algebra $H$ coincides with DAHA of type $A_1$.
For this algebra, the classification of the finite dimensional
representations has been done in \cite{CO}\footnote{Let us remark
that the classification from
 \cite{CO} has a small gap fixed in the book \cite{C}.}.
We classify the finite dimensional representations of the algebra
$H(k_0,k_1,u_0,u_1;q)$ (see  theorem~\ref{main} in the main body of
the paper) and give a description of these representations in terms
of the polynomial representations of $H$. In particular, we
generalize the result from \cite{C}. The proof of the main result of the
paper hinges on  results of Crawley-Boevey on  Deligne-Simpson problem \cite{CB}.

{\bf Acknowledgment} The authors are grateful to Prof. Rogers
for the excellent reseach enviroment  provided by the SPUR-2003 program
during which  this paper was started. We are also grateful to Prof. Etingof
for the useful discussions and valuable references. The work of the
first author was partially supported by the NSF grant DMS-9988796
and the CRDF grant RM1-2545-MO-03.

\section{The double affine Hecke algebra of rank 1}
Throughout this paper, we will work over the field of complex
numbers. Let us consider $q^{1/2}$ a fixed complex number.
\begin{defi}Let $k_0,k_1,u_0,u_1,q\in\CC^*$.
The double affine Hecke algebra $H=H(k_0,k_1,u_0,u_1;q)$ of rank
$1$ of type $CC_1^{\vee}$ is generated by the elements $T_0$, $T_1$,
$T_0^{\vee}$,
$T_1^{\vee}$ with the relations:
\begin{gather*}\label{k0}
(T_0-k_0)(T_0+k_0^{-1})=0,\\
(T_1-k_1)(T_1+k_1^{-1})=0,\label{k1}\\
(T_0^{\vee}-u_0)(T^{\vee}_0+u_0^{-1})=0,\label{u0}\\
(T_1^{\vee}-u_1)(T_1^{\vee}+u_1^{-1})=0,\label{u1}\\
T_1^{\vee}T_1T_0T_0^{\vee}=q^{-1/2}.\label{VVVV}
\end{gather*}
\end{defi}

It is convenient to introduce the following notation:
$$ \underline{t}=(t_1,t_2,t_3,t_4):=(k_0,k_1,u_0,u_1).$$

\subsection{Polynomial representations of $H$}
We give a  description of the finite dimensional irreducible
representations of $H$ in terms of the polynomial representations
$\mathcal{P}^{\pm 1,\pm 1}$ and $\bar{\mathcal{P}}^{\pm 1,\pm 1}$.
Let us describe these representations and their basic properties.

    For each of the subalgebras $H_{aff}=\langle T_1,T_0 \rangle$ and
$H_{aff}^{\vee}=\langle T_1^{\vee}, T_0^{\vee}\rangle$ of $H$, we
define the one-dimensional representations $\chi(\epsilon_0,
\epsilon_1)$  by
$$T_i\rightarrow \epsilon_i k_i^{\epsilon_i}$$ for
$\epsilon_i=\pm 1$, $i=0,1$ and
$\chi^{\vee}(\delta_0, \delta_1)$  by
$$T_i^{\vee}\rightarrow \delta_i k_i^{\delta_i}$$ for
$\delta_i=\pm 1$ and $i=0,1$. We define the following polynomial
representations of the algebra $H$:
\begin{gather*}
\mathcal{P}^{\epsilon_0,\epsilon_1}=H \otimes_{H_{aff}}
\chi(\epsilon_1,\epsilon_0), \\
\bar{\mathcal{P}}^{\delta_0,\delta_1}=H \otimes_{H^{\vee}_{aff}}
\chi^{\vee}(\delta_1, \delta_0).
\end{gather*}

In \cite{NS}, the action of $T_i$ and $T_i^{\vee}$ on the
representation $\mathcal{P}^{1,1}$ is presented. The other cases are
slight modifications of this one. The action of $T_0$, $T_1$ on the
vectors from $\mathcal{P}^{\epsilon_0,\epsilon_1}\simeq
\mathbb{C}[z^{\pm 1}]$ is given by the operators
\begin{gather*}
T_0=\epsilon_0k_0^{\epsilon_0}s_0+
\frac{(k_0-k_0^{-1})+(u_0-
u_0^{-1})q^{1/2}z^{-1}}{(1-qz^{-2})},\\
T_1=\epsilon_1 k_1^{\epsilon_1}s_1+ \frac{(k_1-k_1^{-1})+(u_1-
u_1^{-1})z}{(1-z^{2})}
\end{gather*}
where $s_0(z^n)=q^{n}z^{-n}$, $s_1(z^n)=z^{-n}$ and the element
$q^{1/2}T_0T_0^{\vee}$ acts by multiplication by $z$.

Similarly, the elements $T_0^{\vee}$ and $T_1^{\vee}$ act on
$\bar{\mathcal{P}}^{\delta_0,\delta_1}\simeq \mathbb{C}[z^{\pm 1}]$
by the operators
\begin{gather*}
T^{\vee}_0=\delta_0 u_0^{\delta_0}s_0+
\frac{(u_0-u_0^{-1})+(k_0-
k_0^{-1})q^{1/2}z^{-1}}{(1-qz^{-2})},\\
T^{\vee}_1=\delta_1 u_1^{\delta_1}s_1+
\frac{(u_1-u_1^{-1})+(k_1-
k_1^{-1})z}{(1-z^{2})},
\end{gather*}
and the element $(T_1^{\vee}T_1)^{-1}$ acts by multiplication by $z$.

\section{Classification of the  finite dimensional representations}

\subsection{The root system $D^{(1)}_4$}
The classification uses the affine root system $D^{(1)}_4$  \cite{K}. 
The Dynkin graph for $D^{(1)}_4$ is a
star with four legs of length $1$. The positive roots correspond to
labelings of the vertices of the graph by positive integers with the
special condition explained, for example, in \cite{CB}.

Let $a_1$, $a_2$, $a_3$, $a_4$ be the labels on
the nodes of the four legs, and $a_0$ the label on the central
node of the graph. Thus, we can think of the positive root $\alpha$
as a vector $(a_0,a_1,a_2,a_3,a_4)$ in
$\mathbb{N}^5$.

The positive roots are of two kinds, real and imaginary. We consider
only the strict roots, namely those whose labels of the nodes at the
legs are not smaller than the label of the central node. The
imaginary roots are of the form $(n,n,n,n,2n)=n\Delta$, where $n$ is
a natural number. The real roots are of two types. The strict roots
of the first type are of the form 
$$
r_{i,\epsilon,n}=n\Delta+\epsilon
e_i,$$ 
where $i=0,1$, $\epsilon=\pm 1$, $n>0$ and $e_0=(0,1,0,0,0)$
$e_1=(0,0,1,0,0)$,  or $$r^{\vee}_{i,\delta,n}=n\Delta+\delta f_i,$$
where $i=0,1$, $\delta=\pm 1$, $n>0$ and $f_0=(0,0,0,1,0)$,
$f_1=(0,0,0,0,1)$. The strict roots of the second type are of the
form $$r_{\underline{\epsilon},\underline{\delta},n}=c+n\Delta+
\sum_{i=0}^1\frac{1-\epsilon_i}{2} e_i+
\sum_{i=0}^1\frac{1-\delta_i}{2} f_i,$$ where
$\delta_i,\epsilon_i=\pm 1$, $n>0$ and $c=(1,0,0,0,0)$.

\subsection{Deligne-Simpson problem}
It is easy to see that problem of classification of the
representation of $H$ is closely related  to conjecture~1.4 from the
paper \cite{CB} about the Delign-Simpson problem (they are almost
equivalent). 

Deligne-Simpson problem poses the question of describing the set
of solutions of the system of equations:
\begin{gather*}
A_1\cdot\dots\cdot A_k=1,\\
A_i\in  C_i, \quad i=1,\dots,k,
\end{gather*}
where $C_i\subset GL(n,\CC)$ are the conjugacy classes defined by the equations:
$$\alpha_{ij}=\mathrm{rank}( (A_i-\xi_{i1})\cdot\dots\cdot (A_i-\xi_{ij})),\quad j=1,\dots, w_i,$$
with $\alpha_{i,w_i}=0$ and $\alpha_{i,0}=n$.

Let $\Gamma_w$ be the star-shaped graph with $k$ legs of length $w_1,\dots,w_k.$
The conjecture~1.4 from \cite{CB} claims that there exists an irreducible solution
of the Deligne-Simpson problem if and only if $\alpha$ is a positive root of the root system associated to
$\Gamma_w$ (see \cite{CB} for the definition) and the following equation is satisfied:
\begin{gather*}
\xi^{[\alpha]}=1,\quad \xi^{[\alpha]}:=\prod_{i=1}^k\prod_{j=1}^{w_i}\xi_{ij}^{\alpha_{i,j-1}-\alpha_{ij}}, 
\end{gather*}
together with some with some nonresonance conditions (see \cite{CB} for more details).

If $V$ is an irreducible representation of $H(\underline{t},q)$ then the matrices:
$$ A_1=q^{1/2}T_0,\quad A_2=T_0^\vee,\quad A_3=T_1,\quad A_4=T_1^\vee,$$
is a solution of the Deligne-Simplson problem with $w=(2,2,2,2)$ and the parameters
\begin{gather*}
\xi_{11}=k_0q^{1/2},\quad \xi_{12}=-k_0^{-1}q^{1/2},\\
\xi_{21}=u_0,\quad \xi_{22}=-u_0^{-1},\\
\xi_{31}=u_1,\quad \xi_{32}=-u_1^{-1},\\
\xi_{41}=k_1,\quad \xi_{42}=-k_1^{-1}.
\end{gather*}

\subsection{}To state the main theorem of the paper we need to introduce the
 locally closed subset $\Sigma_\alpha$ ($\alpha$ is a real root of
$D_4^{(1)}$) inside the parameter space $\mathbb{C}^4_{\underline{t}}$:
\begin{itemize}
\item
If $\alpha=r_{\underline{\epsilon},\underline{\delta},n}$
for some $\epsilon_i=\pm 1$ and $\delta_i=\pm 1$, $i=0,1$
then $\underline{t}\in \Sigma_\alpha$ iff
$$\epsilon_1 k_1^{\epsilon_1}\epsilon_0
k_0^{\epsilon_0}\delta_1 u_1^{\delta_1}\delta_0 u_0^{\delta_0}=q^{-1/2-n}$$
and
$k_i^{2\epsilon_i}\neq -q^{m_i}$, $u_i^{2\epsilon_i}\neq -q^{m_i}$ for
all $i=0,1$, and $(1+\epsilon_i)/2\le m\le n-1$.

\item
If  $\alpha=r_{i,\epsilon,n}$
for some $i=0,1$ and some $\epsilon=\pm 1$
then $\underline{t}\in
\Sigma_\alpha$  iff
 $$k_i^{2\epsilon_i}=-q^{n}$$  and
$\epsilon_1 k_1^{\epsilon_1}\epsilon_0 k_0^{\epsilon_0}\delta_1
u_1^{\delta_1}\delta_0 u_0^{\delta_0} \neq q^{-1/2-m}$
for all $m=0,1,...,n-1$

\item
If  $\alpha=r^{\vee}_{i,\delta,n}$ for some $i=0,1$ and some $\delta=\pm 1$
then $\underline{t}\in
\Sigma_\alpha$  iff
 $$k_i^{2\epsilon_i}=-q^{n}$$  and
$\epsilon_1 k_1^{\epsilon_1}\epsilon_0 k_0^{\epsilon_0}\delta_1
u_1^{\delta_1}\delta_0 u_0^{\delta_0} \neq q^{-1/2-m}$
for all $m=0,1,...,n-1$
\end{itemize}

Given a root $\alpha=(a_0,a_1,a_2,a_3,a_4)$ we can set $\alpha_{i1}=a_1$ and $\alpha_{i0}=a_0$ then the
the equation for the closure $\bar{\Sigma}_\alpha$ is of the form $\xi^{[\alpha]}=1$.

\subsection{Classification}

Let us  define the polynomials
\begin{gather*}
E_n(a;z)=z^{-n}\prod_{i=-n}^{-1}(z-q^i/a)\prod_{i=0}^n(z-aq^i)\\
E_{-n}(a;z)=z^{-n}\prod_{i=-n}^{-1}(z-q^i/a)\prod_{i=0}^{n-1}(z-aq^i).
\end{gather*}

Suppose that $V$ is a finite dimensional representation of
$H(\underline{t},q)$. Then we define $\dim(V)=
\vec{d}=(d_0,d_1,d_2,d_3,d_4)$ where $d_0=\dim(V)$,
$d_1=\dim(Im(T_0-k_0))$, $d_2=\dim(Im(T_1-k_1))$,
$d_3=\dim(Im(T_0^\vee-u_0))$, $d_4=\dim(Im(T_1^\vee-u_1))$.

\begin{thm}\label{main} Suppose that $q$ is not a root of $1$. Then
\begin{itemize}
\item The algebra $H=H(\underline{t};q)$ has a finite dimensional
irreducible representation of dimension $\vec{d}$ iff $\alpha=
d_0c+d_1e_0+d_2e_1+d_3f_0+d_4f_1$ is a strict real  root for the
$D_4^{(1)}$ root system and $\underline{t}\in\Sigma_\alpha$.
\item The algebra $H$ could have only one
 irreducible representation
 of dimension $\alpha$. This representation $T_\alpha$ has an
explicit description as a quotient  of the polynomial representation:
\begin{align*}
V_\alpha=&\mathcal{P}^{\epsilon_0,\epsilon_1}/(E_n(\epsilon_0
k_0^{\epsilon_0} \delta_0 u_0^{\delta_0};z))\quad &\mbox{ when } 
\quad \alpha=r_{\underline{\epsilon},\underline{\delta},n}\\
V_\alpha=&\mathcal{P}^{\epsilon,1}/ (E_{-n}(q^{1/2}\epsilon
k_0^{\epsilon}u_0;z)),
\quad&\mbox{ when }\quad \alpha=r_{0,\epsilon,n},\\
V_\alpha=&\mathcal{P}^{1,\epsilon}/ (E_{-n}(q^{1/2}\epsilon
k_1^{\epsilon}u_1;z)),\quad &\mbox{ when }\quad \alpha=r_{1,\epsilon,n},\\
V_\alpha=&\bar{\mathcal{P}}^{\delta,1}/ (E_{-n}(\delta
u_0^{\delta}k_0;z)),
\quad&\mbox{ when }\quad \alpha=r_{0,\delta,n},\\
V_\alpha=&\bar{\mathcal{P}}^{1,\delta}/ (E_{-n}(\delta
u_1^{\delta}k_0;z)),\quad&\mbox{ when }\quad \alpha=r_{1,\delta,n}. \\
\end{align*}
\end{itemize}

\end{thm}

\begin{proof}
In the next section we prove that all finite dimensional
representations of $H$ are rigid (see Lemma~\ref{rig}). Hence the
first part of the theorem is equivalent to theorem~1.5 from
\cite{CB}.

Let us give the proof for the first case of the second half of the theorem.
If we have a representation $V$ of  dimension
$\alpha=r_{\underline{\epsilon},\underline{\delta},n}$ then from the definition
of  dimension we get that there  exist vectors $w,v\in V$ such that
\begin{gather*}
T_0w=\epsilon_0 k_0^{\epsilon_0}w,\quad T_1w=\epsilon_1 k_1^{\epsilon_1} w,\\
T_0^{\vee}v=\delta_0u_0^{\delta_0} v,\quad T_0v=\epsilon_0k_0^{\epsilon_0}v.
\end{gather*}

Indeed, if we know that $k_i\ne -k_i^{-1}$, $u_i\ne -u_{i}^{-1}$,
$i=0,1$  then $T_i, T_i^{\vee}$ are semisimple and the statement
follows from the first part. If say $k_i=-k_i^{-1}$ for some $i=0,1$
then the first part implies $\underline{t}\in \Sigma_{\alpha}$.
Hence $\epsilon_i=1$ and $\dim(ker (T_i-\epsilon_ik_i)) =n+1$.

This implies that the map $\mathcal{P}_{\underline{\epsilon}}\ni 1\mapsto v\in
V$ extends to an
$H$-homomorphism $\varphi$: $\mathcal{P}^{\underline{\epsilon}}\to
V$. The kernel of $\varphi$ is the ideal $(E)$, where
$E\in\mathbb{C}[z^{\pm 1}]=\mathcal{P}_{\underline{\epsilon}}$ is a Laurent
polynomial.

To identify $E$ we need to find the spectrum of $z=q^{1/2}T_0T_0^{\vee}$ on
$V$. It is easy because for the same reason as before we have the
$H$-homomorphism $\psi$: $V\to H\otimes_{H'}\chi=\mathcal{P}'$ where
$H'$ is a subalgebra of $H$ generated by $T_0,T_0^{\vee}$:
\begin{equation*} \chi(T_0)=\epsilon_0 k_0^{\epsilon_0},\quad \chi(T_0^{\vee})=
\delta_0 u_0^{\delta_0}.
\end{equation*}

Because of the PBW theorem we have a canonical identification
$\mathcal{P}'\simeq\mathbb{C}[y^{\pm 1}]$ with $y=T_1T_0$.
Analogously to the polynomial representations
$\mathcal{P}^{\underline{\epsilon}}$ (see for example \cite{NS}
section~3.5) we can introduce a complete ordering on the monomials:
$$1\prec y^{-1}\prec y\prec y^{-2}\prec y^2\prec \dots.$$
The action of the operator $z=q^{1/2}T_0T_0^{\vee}$ is upper
triangular with respect to this ordering:
$$z(y^i)=\rho_iy^i+\mbox{lower order terms},$$
where $\rho_i=\epsilon_0\delta_0 k_0^{\epsilon_0}u_0^{\delta_0}q^{1/2+i}$
for $i\ge 0$ and
      $\rho_i=\epsilon_0\delta_0 k_0^{-\epsilon_0}u_0^{-\delta_0}q^{1/2+i}$ for $i<0$.

Thus we have the identification $V\simeq \mathcal{P}'/(E')$ where $E'\in\mathbb{C}[y^{\pm 1}]$.
We know that $\dim(V)=2n+1$, hence we can choose $E'$ in the form
$$E'=y^n+cy^{-n}+\mbox{lower order terms},$$
with $c\ne 0$. This implies that  the monomials $1,y^{-1},y,^1,\dots,y^{-n},y^n$ project onto a
basis in $\mathcal{P}'/(E')$. The
upper triangularity of $z$ implies that the spectrum of
$z$ is $\rho_{-n},\dots,\rho_n$. This implies $E=E_n(\epsilon_0 k_0^{\epsilon_0}
\delta_0 k_0^{\delta_0};z)$. It can be shown by a direct computation that
the ideal $(E_n(\epsilon_0 k_0^{\epsilon_0}
\delta_0 u_0^{\delta_0};z))$ is invariant under the action of $T_0,T_1$ if $\underline{t}\in
\Sigma_{\alpha}$.

\end{proof}

\subsection{Rigidity of the finite dimensional representations}

In this section we prove
\begin{lem}\label{rig} Suppose that $q$ is not a root of unity and the
algebra $H=H(\underline{t};q)$ has a finite dimensional
irreducible representation $V$. Then this representation is rigid,
i.e. it does not admit any deformations.
\end{lem}

The proof of Lemma~\ref{rig} goes along the lines  of the  calculation from
 theorem 5.7 in \cite{EOR}).

\begin{prop} Let $C_i\subset GL(n)$, $i=1,\dots,4$ be
the conjugacy classes in $GL(n)$ (i.e. subset of matrices with fixed
Jordan normal form) and $\mathcal{C} = C_1\times C_2 \times
C_3\times C_4$ be the subset of the irreducible 4-tuples. Let
$\Delta\subset \mathcal{C}$ be the subset of 4-tuples $(a,b,c,d)$ of
matrices with the property $abcd=q^{1/2}$. If $\det(C_1)\det(C_2)
\det(C_3)\det(C_4)=q^{n/2}$ then
$$\dim\Delta=\dim C_1+
\dim C_2+\dim C_3+\dim C_4-n^2+1,$$ when RHS is nonegative and
$\Delta=\emptyset$ if RHS is negative.
\end{prop}
\begin{proof}
Obviously we have the equality
$$ \dim\mathcal{C}=\dim C_1+
\dim C_2+\dim C_3+\dim C_4,$$ and inequality
$$\dim \Delta\ge \dim\mathcal{C}-(n^2-1).$$ To prove the opposite
inequality we need to show that the value $q^{1/2}d^{-1}$ is
regular for the map $\nu: C_1\times C_2\times C_3\to GL(n)$ given
by $\nu(a,b,c)=abc$. Let us recall that in our case the value
$q^{1/2}d^{-1}$ is regular if and only if the dimension of the
image of the linear map $d\nu$ is equal to $n^2-1$. Let us explain
why this is true.

On the space of matrices there is a natural  nondegenerate pairing
$x,y\mapsto tr(xy)$. Let us use this pairing to describe the
orthogonal complement to the image of the differential of the map
$\nu$.
 It is easy to see that $t$ belongs to the orthogonal complement
 if and only if
it satisfies equations \begin{equation} \label{orth}[a,bct]=0,
\quad [cta,b]=0,\quad [tab,c]=0.\end{equation}  Indeed let
$u=[x,a]$, $v=[y,b]$, $w=[z,c]$ be the tangent vectors to
$C_1,C_2,C_3$ at the points $a,b,c$. Then $t$ satisfies the
equation $tr(t([x,a]bc+a[y,b]c+ab[z,c]))$ for any $x,y,z$. Now the
cyclic invariance of the trace implies the equations (\ref{orth}).

Set $r=bcta$. Then $[a,r]=[b,r]=[c,r]=0$. Since $a,b,c$ are
irreducible, $r$ is a scalar, and $t=\lambda a^{-1}b^{-1}c^{-1}$
for some $\lambda\in\mathbb{C}$. That is, the value $q^{1/2}d^{-1}$
is regular.
\end{proof}

\begin{proof}[Proof of Lemma~\ref{rig}]
Suppose $V$ is a representation of the algebra
$H(\underline{t};q)$ of dimension $n$. Let $C_1$, $C_2$, $C_3$,
$C_4$ be the conjugacy classes of the elements
$T_0,T_0^\vee,T_1^\vee,T_1\in End(V)$.  The quadratic relations
for the elements $T_0,T_0^\vee,T_1^\vee,T_1$ give the restrictions
on the possible type of the Jordan normal form of the elements from
$C_i$, $i=1,\dots,4$. Namely there are two cases.

The first case is when $t_i^2\ne -1$. In this case all matrices
from $C_i$ are diagonalizable. That is, there is a  number $1\le
d_i\le n$ such that for all matrices $X\in C_i$ we have $\dim
Ker(X-t_i)=d_i$ and $\dim Ker(X+t_i^{-1})=n-d_i$. In this case
$\dim C_i=2d_i(n-d_i)$.

The second case is when $t_i^2=-1$. In this case there is a number
$1\le d_i\le n$ such that $\dim Ker(X-t_i)=d_i$ for all $X\in C_i$.
That is, the Jordan normal form of $X$ has $n-d_i$ Jordan block of
size $2$. In this case we have, again, $\dim C_i=2d_i(n-d_i)$.

Now let us remark that the subset $\Delta$ from the previous lemma is
the subset of $n$-dimensional irreducible representations of $H(\underline{t};q)$ with
fixed Jordan normal forms of elements $T_0,T_0^\vee,T_1^\vee,T_1$.
The action of $PGL(n)$ on the space of $n$-dimensional irreducible representations is
free, hence the space of equivalence classes of irreducible representations
$\Delta/PGL(n)$ has the dimension:
$$D=D(d_1,d_2,d_3,d_4)=2(1-n^2+\sum_{i=1}^4 d_i(n-d_i)),$$
where $d_1=\dim(Ker(T_0-t_1))$, $d_2=\dim(Ker(T_0^\vee-t_2))$,
$d_3=\dim(Ker(T_1^\vee-t_3))$, $d_4=\dim(Ker(T_1-t_4)).$

It is easy to see that if $n=2N+1$ then $D\ge 0$ if and only if
$d_i=N+\delta_i$, $\delta_i=0,1$.  In this case $D=0$ and the
calculation of the determinants gives $\prod_{i=1}^4(\epsilon_i
t_i^{\epsilon_i})=q^{1/2+N}$, $\epsilon_i=(2d_i-1)/2$,
$i=1,\dots,4$. That is, this irreducible represntation is rigid.

If $n=2N$, there are two possibilities for $D$ to be positive. The
first case is when $d_i=N$. In this case we have $D=2$ and
calculation of the determinants gives $q^N=1$. The second case is
when there  exists $1\le i\le 4$ and $\epsilon=\pm 1$ such that
$d_i=N+\epsilon$ and $d_j=N$ for $j\ne i$. In this case $D=0$ and
$t_i^{2\epsilon}=-q^N$.
\end{proof}


\begin{thebibliography}{99}
\bibitem{C1} I.Cherednik, \emph{Double affine Hecke algebras, KZ
equations and Macdonald operators}, IMRN (1992), no. 9, 171--180.

\bibitem{C2} I.Cherednik, \emph{Macdonald's evaluation conjectures and
difference Fourier transform}, Invent. Math. vol. 122 (1995), no.1, 191--216.

\bibitem{EG}P.Etingof, V.Ginzburg, \emph{ Symplectic reflection algebras,
Calogero-Moser space, and deformed Harish-Chandra homomorphism},
Invent. Math., vol.147 (2002), no. 2, 243--348.

\bibitem{Sa} S.Sahi, \emph{ Nonsymmetric Koornwinder polynomials and duality}, Ann. of Math. (2),
150, (1999), no.1, 267-282.


\bibitem{NS} M.Noumi, J.V.Stokman, \emph{Askey-Wilson polynomials: an
affine Hecke algebraic approach}, math.QA/0001033 (2000).

\bibitem{O} A.Oblomkov, \emph{Double affine Hecke algebras and
and affine cubic surfaces}, IMRN, no. 18, 877-912.

\bibitem{CO} I.Cherednik, V.Ostrik, \emph{From double affine Hecke
algebra to Fourier transform}, math.QA/0111130 (2001).

\bibitem{C} I.Cherednik, \emph{Double affine Hecke algebras}, Cambridge University Press, 2004.

\bibitem{CB}W.Crawley-Boevey, \emph{Indecomposable parabolic bundles
and the existence of matrices in prescribed conjugacy class closures
with product equal to the identity}, preprint math.AG/0307246.

\bibitem{EOR} P. Etingof, A. Oblomkov, E. Rains, \emph{
Generalized double affine Hecke algebras of rank 1
and quantized Del Pezzo surfaces}, in preprint QA/0406480.

\bibitem{K} V. Kac, \emph{Infinite Dimensional Lie Algebras - an introduction},
Birkh\"{a}user, 1984, 38-44.


\end{thebibliography}
\end{document}